\input amstex
\documentstyle{amsppt}

\loadeufb
\loadeusb
\loadeufm
\loadeurb
\loadeusm

\magnification =\magstep 1
\refstyle{A}
\NoRunningHeads

\topmatter
\title  On ampleness of canonical bundle  
\endtitle
\author Robert  Treger \endauthor
\address Princeton, NJ 08540  \endaddress
\email roberttreger117{\@}gmail.com \endemail
\endtopmatter

\document

In the present  note we extend the results of \cite{Tr1}. We assume the reader is familiar with \cite{Tr1}. The main new observation of the note is the lemma below.

Let $\phi: X\hookrightarrow \bold P^r$ be a projective manifold of dimension $n\geq 1$. Let $U_X$ denote 
the universal covering of $X$.  We assume $U_X$ is equipped with an arbitrary real analytic
$\pi_1(X)$-invariant Kahler metric. Until Theorem 2, we assume $\pi_1(X)$ is nonamenable.

 Let $R$ be a Riemannian manifold which is a Galois covering of a compact manifold $N$ with nonamenable Galois group. 
In the fundamental paper \cite{LS, Theorems 3, 3$'$}, Lyons and Sullivan proved  $R$ admits a non-constant bounded harmonic function.
Employing their theorem, Toledo proved 
that the space of bounded harmonic functions on $R$
 is infinite dimensional \cite{To1} (see a brief review in \cite{Tr1, (2.6)}). 

Assuming $\pi_1(X)$ is nonamenable,
let $V^b$ denote the vector space generated by all
bounded positive pluriharmonic functions on $U_X$ with the sup norm. 
The fundamental group $\pi_1(X)$ acts on $V^b$ by isometries.  We get 
a representation
$$
\rho: \pi_1(X) \longrightarrow Isom( V^b).
$$
Let $\Xi$ denote the kernel of $\rho$.  Let $U:=U_X/\Xi$  denote a Galois covering of $X$ with the Galois group denoted by $\Gamma$.
Our Kahler metric on $U_X$ induces the $\Gamma$-invariant real analytic Kahler 
metric on U. Let $V^b_U$ denote the vector space generated by all {\it  bounded}\/ 
positive pluriharmonic functions on $U$ with the sup norm. As before, the spaces $V^b_U$ and $V^b$ are naturally isomorphic {\it infinite-dimensional}\/ vector spaces.

Let $PHar(U_X)$  
be the vector space generated by {\it all}\/ positive pluriharmonic functions  
on $U_X$. 
We will integrate   pluriharmonic functions  
with respect to the measure
$$
dv:=p_{U_X}(s,x,\bold Q) d\mu = p_{U_X}(x) d\mu,
$$
where $\bold Q \in U_X$ is a fixed point, 
$d\mu$ is the corresponding  Riemannian measure, and 
$p_{U_X}(x):= p_{U_X}(s,x,\bold Q)$ is the corresponding heat kernel. 
We obtain a pre-Hilbert space generated by the {\it bounded}\/ positive pluriharmonic 
square integrable functions on $U_X$ (compare \cite{Tr1, Sect.\;2.4 and Sect.\;4}).

The latter pre-Hilbert space  has a completion in the Hilbert space 
$H$ generated by all positive pluriharmonic  $L_2(dv)$ functions:
$$
H := \biggl\{h \in PHar(U_X)\;\biggl |\; \parallel h\parallel^2_H:=  \int_{U_X}
|h (x)|^2 dv =\int_{U_X}|h (x)|^2 p_{U_X}(x) d\mu < \infty   \biggl\}.
$$ 
Let $H^b \subseteq H$ be the Hilbert subspace generated by $V^b$.
These Hilbert spaces are separable infinite-dimensional Hilbert spaces with reproducing kernels.

Similarly, we consider the vector space $PHar(U) $, the corresponding
heat kernel $p_U(s.x,\bold Q)$, the measure $dv_U$ in place of
$dv$,  the Hilbert space $H_U$ in place of
$H$, and the Hilbert subspace $H_U^b \subseteq H_U$ generated by $V^b_U$.

Given $u\in V_U^b$, there exists a holomorphic $L_2(dv_{U})$ function $f=u+
\sqrt{-1}\tilde u$ on $U$. Namely, 
we set $\tilde u(x)= \sqrt{-1}\int_{x_0}^x (\bar\partial u - \partial u)$,
where $x_0 \in U$ is a fixed point and $x\in U$ is a variable point
(see, e.g., \cite{FG, Chap.\;6.1, p.\;318} and \cite{Tr1, Sections 2.5, 4.1, 4.2}).

Similarly to \cite{Kl1, Prop.\;2.12}, $X$ is said to be $\Gamma$-large if $U$ contains no compact positive dimensional analytic subsets.

\proclaim{Theorem 1 (Uniformization II)} We keep the above notation. If  
$\pi_1(X)$ is non-amenable, $\Gamma$-large  and non-residually finite then the canonical bundle $\eusm K_X$ is ample and $U$ is a bounded Stein domain in $\bold C^n$.
It follows $U_X$ is Stein as well.
\endproclaim

\demo{Proof} Let  $\Cal L$ by a very ample line bundle defining $\phi$.
As in \cite{Tr1}, we can construct the real analytic Kahler metrics $\Lambda_{U,\Cal L}$, $\Sigma_{U,\Cal L} $ and $\beta_U$ on $U$ in place of $\Lambda_{\Cal L}$, $\Sigma_{\Cal L} $ and $\beta$ on $U_X$.
Then we show show that $U$ is a bounded domain in $\bold C^n$, provided we can establish that $\Gamma$ is {\it residually finite}\/ (see Lemma below). 

We derive that  $U$  is Stein by Siegel's theorem; see a discussion of Siegel's theorem and  references in \cite{Tr1, (A.0)}.
Because  $U_X$ is a covering of $U$, $U_X$ is Stein as well.

The Prolongation Lemma (see \cite{Tr1, Lemma A in Appendix})  
is valid on $U$. Indeed, the diastasic potential on $U$ is induced by
the  diastasic potential in the target Fubini space. Furthermore, 
the  Bochner canonical coordinates in $U$ are  holomorphic functions on the
whole $U$.

In the proof that $U$ is a bounded domain in $\bold C^n$, we proceed as in \cite{Tr1,
(5.3)} which rely on the fundamental paper by Calabi \cite{C, Theorem 7 on p.\;15,
Proposition 7 on p.\;14, Theorem 6 on p.\;13, and Theorem 12 and its proof on pp.\;20-21}.
We observe that the natural image of $U$ in the Fubini space $\bold F_\bold C(\infty,1)
$ 
does not intersect the corresponding antipolar hyperplane 
because,  in our case,
the diastasic potential of $U\hookrightarrow \bold F_\bold C(\infty,1)$ is a {\it 
 function}\/  on $U$.

This completes the prove of the theorem provided $\Gamma$ is residually finite.
\enddemo

\remark{Remark 1} Recall that Toledo constructed an example of projective manifold with nonamenable and non-residually finite fundamental group \cite{To2}.

If  $M$ is an arbitrary compact real analytic Kahler manifold with  generically large and
nonamenable fundamental group then $M$  is projective by \cite{Tr3} (there we have proved that  $M$ is Moishezon hence projective). Recall that if the corresponding fundamental group is residually finite then  the well-known
H.\;Wu conjecture about Kahler manifolds with negative sectional curvature is valid \cite{Tr2}.
\endremark

\proclaim{Corollary} With notation and assumptions of the theorem,  
$U_X$ is not a bounded domain in $\bold C^n$.
\endproclaim

\demo{Proof} Suppose $U_X$ is a bounded domain in $\bold C^n$. 
Then the functions of $V^b$ separate points on $U_X$. Since $V^b$ and $V_U^b$ are naturally isomorphic and $U_X\rightarrow U$ is a nontrivial covering, we get a contradiction.
\enddemo

To complete the proof of the theorem and its corollary, we need the following lemma which is a generalization of a classical theorem of Maltsev about 
finitely generated subgroups of $GL(m)$ for
$m < \infty$ (see, e.g., \cite{Z, Chap.\;1.2}).

\proclaim {Lemma} Let  $V$ be a  vector  space over  $\bold C$  with a countable base. Let $\eusb G(\bold C)$ denote the group of all $\bold C$-linear automorphisms of $V.$ Let  $H \subset \eusb G(\bold C)$ be a finitely generated subgroup. Then $H $ is residually finite.
\endproclaim

\demo {Proof} 
We will assume $\dim V=\infty$.  The case $\dim V < \infty$ was treated by Maltsev.
We fix a base in $V$. 
Each $g\in \eusb G(\bold C)$ is given by an $(\infty \times \infty)$ matrix with entries in $\bold C$. 
Given an  $(\infty\times \infty)$ matrix $M$, let  $M_i\; (1\leq i<\infty)$ 
denote the matrix whose all entries outside the upper left  $(i\times i)$ corner
 block of $M$ are replaced by zeros.

Roughly speaking, all $(\infty \times \infty)$ matrices $(x_{jk})\; (1\leq j,k <\infty)$ with indeterminate entries have a structure
of  an affine ind-algebraic variety  $\eufm M$ (an affine infinite-dimensional variety 
in the sense of Shafarevich \cite{S})
defined by finite-dimensional subvarieties $\eufm M_i$ in the obvious way, where 
$\eufm M_i =\{M_i\}$. 

The  algebra  $\bold C[\eufm M]$
of regular functions on $\eufm M $ has a structure of   topological algebra.
Let  $\bold C[\eufm M\times \eufm M]=\bold C[\eufm M] \hat \otimes_\bold C \bold C[\eufm M]$ be the algebra of regular functions on $ \eufm M \times \eufm M$. 
The regular functions on $\eufm M_i$
form an algebra and a coalgebra with the standard comultiplication
$$
\Delta_i : \bold C[\eufm M_i] \rightarrow \bold C[\eufm M_i  \times  \eufm M_i],
\; (\Delta_i f_i)(g'_i,g''_i) =f_i(g'_i g''_i)\quad ( f_i\in \bold C[\eufm M_i];\;
g'_i, g''_i \in \eufm M_i).
$$
We consider each  $\eufm M_i$  with its Zariski  topology, and $\eufm M$ is equipped with the topology of inductive limit \cite{S, Sect\;1}. In fact,
$\eufm M_i$ is a finite-dimensional algebraic semigroup variety. 
Since $\eufm M=\varinjlim \eufm M_i$, we get a natural multiplication map
$\eufm M \times \eufm M \rightarrow \eufm M$ associated with a 
natural map of topological algebras
$
\Delta_\eufm M : \bold C[\eufm M] \rightarrow \bold C[\eufm M  \times\eufm M], 
$ which is the key point of the lemma. 
Let $\eufm M(\bold C)$ be the $\bold C$-points of $\eufm M$. 

A set $G$ which is a group and an ind-algebraic variety
is said to be an ind-algebraic group if the inversion map $G\rightarrow G$ and the  
multiplication map $G\times G \rightarrow G$ are morphisms of ind-algebraic varieties (\cite{S, Sect.\;1}, \cite{AT, Sect.\;3.2}).

Let $\eurb e$ be the $(\infty\times\infty)$  unit matrix. 
Let $\eufm M'$ and $\eufm M''$ denote two copies of $\eufm M$.
The structure of ind-affine variety on $\eufm M' \times \eufm M''$ is defined by 
finite-dimensional subvarieties $\eufm M_i' \times \eufm M_i'' 
\subset \eufm M' \times \eufm M''$, where $\eufm M'_i$ and $\eufm M''_i$ are two copies of $\eufm M_i$.

The  equation  $\eufm p \circ \eufm q =\eurb e$  
($\eufm p\in\eufm M', \eufm q \in \eufm M''$) 
defines an ind-algebraic group as well as an ind-algebraic subvariety 
$\eufm G \subset \eufm M' \times \eufm M''$ (compare  \cite {S, Sect.\;2};
the equation produces the projective system of algebras of the form
$R_i=\bold C[x'_{jk}; x''_{jk}]/\{\text{relations}\}$ defining $\eufm G$). 
From the projection $\eufm M' \times \eufm M''\rightarrow \eufm M' $, we get the embedding $\eufm G \subset \eufm M $.

Let $\eufm G(\bold C)$  denote the $\bold C$-points of
$\eufm G$. We get a map of groups
$\eta_H: H\rightarrow \eufm G(\bold C)$; $\eta_H$ arises from the map $H\rightarrow (H,H),\/\;  h\mapsto (h,h^{-1})$.
The closure of $\eta_H(H)$ in $\eufm G(\bold C)$, denoted by $\overline{\eta_H(H)}$,
  is a group as in the finite-dimensional case (compare \cite{M, Lemma 1.2.6}).
Therefore $\overline{\eta_H(H)} \subset \eufm G(\bold C)\subset \eufm M(\bold C)$ 
is an ind-algebraic subgroup as in the finite-dimensional case. In fact, the algebra of regular functions on $\overline{\eta_H(H)}$ is a topological Hopf algebra with the standard comultiplication and antipode and the structure maps satisfy the well-known identities. 

Let $\eusb M(\bold C)$ be the affine space of $(\infty\times\infty)$ matrices with entries in $\bold C$ with its Zariski topology. 
Let $I: \eufb M(\bold C) \rightarrow \eusb M(\bold C)$ be the natural continuous map. 

Let $\eusb H(\bold C)$ be the closure of $H$ in $I(\eufm G(\bold C))$.
Clearly, $I$ is an isomorphism of the group $H\subset\eufm M(\bold C)$ onto the group $H\subset I(\eufm G(\bold C))$. Also, $I$  
maps $\overline{\eta_H(H)}$ onto  $\eusb H(\bold C)$. 

Indeed, given $ h\in \eusb H(\bold C)\backslash H$, we take a general  curvilinear section $C\; (\subset \eusb M(\bold C))$ though $ h$. Then $I^{-1}(C\cap \eusb H(\bold C))$ will be closed in  $\overline{\eta_H(H)}$
and $I^{-1}( h) $ will be a unique point in $ \overline{\eta_H(H)} $, i.e., $I$ maps
$\overline{\eta_H(H)}$ one-to-one onto $\eusb H(\bold C)$.
Further, $\eusb H(\bold C)$ is an open subset of its closure in $\eusb M(\bold C)$
(Generalized Chevalley theorem).

Let $\Cal A$ and $\Cal B$ denote the commutative (reduced) $\bold C$-algebras of regular functions on $\eusb H(\bold C)$ and $\eusb H(\bold C)\times \eusb H(\bold C)$, respectfully. 
Our aim is to show that $\eusb H(\bold C)$ has a structure of an {\it affine group}\/ \cite{AT, Chap.\;3.2}. 
We will show $\Cal A$ has a natural structure of a discrete commutative Hopf algebra.

Let $A(H):=\Cal A|_H$ and $B(H\times H): =\Cal B|_{H\times H}$ denote
 $\bold C$-algebras of the corresponding $\bold C$-maps $H\rightarrow \bold C$
and $B(H\times H)  \rightarrow \bold C$, respectfully 
($ H\times H$ is a subset of $\eusb H(\bold C)\times \eusb H(\bold C)$).
We observe that
$
\Cal B =  \Cal A\otimes_\bold C \Cal A$ hence  $B(H\times H) = A(H)\otimes_\bold C A(H).
$
It follows  the existence of 
 a discrete Hopf algebra structure on  $A(H)$ with the co-identity $\epsilon_H$, the
comultiplication $\Delta_H$ and the antipod $S_H$ (see, e.g., \cite{M, Chap.\;3}),
$$ 
\epsilon_H f_H:=f_H(1),\quad \Delta_H:  A(H) \longrightarrow  A(H\times H), \; (\Delta_H f_H)(h',h''):=
f_H(h'h''),
$$
$$
S_H: A(H) \longrightarrow  A(H), \; (S_Hf_H)(h): = f_H(h^{-1}) \quad 
(f_H\in A(H); \; h, h', h'' \in H).
$$

A regular function $f_\Cal A$ on $\eusb H(\bold C) $ in the infinite number of variables
$\{x_{jk}\}$ defines a projective system $\{f_{\Cal A,i}\}$, where each $f_{\Cal A,i}$ 
is obtained from $f_\Cal A$ by letting  the variables $\{x_{jk}\}$ equal to zero when
$j>i$ or $k>i$.

Now, we consider arbitrary $f_\Cal A\in \Cal A$ and $h', h''\in \eusb H(\bold C)$. 
Let $h'=\lim h'_t$ and $h''=\lim h''_t \; (h'_t, h''_t \in H;\; t\in \bold N)$. 
As before, we  assume $\{h'_t, h'\}$ as well as $\{h''_t, h''\}$ are contained in the corresponding curvilinear sections. Clearly,
 $I^{-1}(h'_t ) I^{-1}(h''_t )= I^{-1}(h'_th''_t )   $.
It follows
$\lim I^{-1}(h_t'),\; \lim I^{-1}(h_t'')$ and $ \lim I^{-1}(h_t' h_t'')$ exist.

Furthermore,  $(\lim I^{-1}(h_t))^{-1} = \lim ((I^{-1}(h_t))^{-1})=\lim I^{-1}(h_t^{-1})  $ and the limits exist where $h_t:=h'_t$. The $\lim I^{-1}(h_t')\lim I^{-1}(h_t'')$ makes sense because $\overline{\eta_H( H) }$ is a group.

Let $\hat \Cal A$ be the topological Hopf algebra of regular functions on 
$\overline {\eta_H(H)}$; we have 
$\hat \Cal A = \varprojlim \bold C[\overline{\eta_H( H) }_i ]$.
Let $f_{\hat \Cal A}:= \varprojlim f_{\Cal A,i}$ be the image of $f_ \Cal A $ in $\hat \Cal A$. We get
$$
 (\Delta_{\hat\Cal A} f_{\hat \Cal A})(\lim I^{-1}\!(h_t'), \lim I^{-1}\!(h_t''))\!=\!
 f_{\hat \Cal A}(\lim I^{-1}\!(h_t')\lim I^{-1}\!(h_t''))\!=\!
\lim f_{\hat \Cal A}(I^{-1}\!(h_t'h_t'')). 
$$

Hence, we can define a comultiplication in $\Cal A$ as follows:
$$
(\Delta_\Cal A f_\Cal A)(h', h'') = (\Delta_\Cal A f_\Cal A)(\lim h_t', \lim h_t'') :=
\lim f_{\Cal A}(h_t'h_t'').
$$
The definition is independent of choices of $h'_t$ and  $h''_t$.
Similarly, we define the antipode $S_\Cal A$. Set $\epsilon f_\Cal A:= f_H(1)$.
The structure maps satisfy the well-known identities.

Thus we obtain a  discrete  Hopf $\bold C$-algebra structure on $\Cal A$ 
extending the 
Hopf algebra $A(H)$.
Hence 
$\eusb H(\bold C)$ has a structure of an affine  group associated with  $\Cal A$ (see \cite{AT, Sect.\;3.2}). 
We get $\Cal A =\varinjlim \Cal A_\alpha$ is the inductive limit of finitely 
generated sub-Hopf algebras $\Cal A_\alpha \subset \Cal A$ (see \cite{A, Lemma 3.4.5}). 
Hence $\eusb H(\bold C)=\varprojlim \eusb H_\alpha(\bold C)$  is the projective limit
of (finite-dimensional) affine algebraic groups $\eusb H_\alpha(\bold C) $.

So $ \eusb H(\bold C)$ is a pro-affine algebraic group. The projection of $H$  in each $\eusb H_\alpha(\bold C)$ is residually finite by  Maltsev's theorem. 
Hence $H$ is residually finite.
\enddemo

The next theorem is a higher-dimensional generalization of the Poincar\'e
ampleness theorem ($\dim X=1)$. In the sequel, $\pi_1(X)$ is not necessary nonamenable.

\proclaim{Theorem 2} We keep the above notation. Let $X$ be a projective manifold with residually finite and large fundamental
group. If the genus $g(C)$ of a general curvilinea section $C\subset X$ is at least $2$ then the canonical bundle $\eusm K_X$ is ample.
\endproclaim

\demo{Proof} Let  $\Cal L$ by a very ample line bundle defining $\phi$.
We can construct the real analytic Kahler metrics $\Lambda_\Cal L$ 
(a generalization of Poincar\'e metric) \cite {Tr1, Sect.\;3.3}.
Let $D_{U_X}(\bold Q,p)$ denote the diastasic potential at $\bold Q \in U_X$ of our
Kahler metric. 

We will follow \cite {Tr1} with some corrections.
Because the diastasis is inductive on complex submanifolds, various questions about higher-dimensional manifolds are reduced to the one-dimensional case (a
generalization of the Poincar\'e metric, the proof of Prolongation lemma, the Shafarevich conjecture, etc. \cite{Tr1}). 

We consider the real-valued function 
$\tilde \Phi _\bold Q (z(p),\bar z(p)):=D_{U_X}(\bold Q,p)$ in a small 
neighborhood $\Cal V_\bold Q \subset U_X$.

Let $C$ be a  general curvilinear section of $X$ and $R$ is its inverse image on $U_X$. With the notation from \cite{Tr1, Proposition-Definition 2}, we can construct a real 
analytic  $Gal(R/C)$-invariant Kahler metric $\Lambda_R :=\lim_{t \to \infty}{1\over t} g_{R,t}$.  
It  follows from the  prolongation over  $U_X$ \cite{Tr1, Appendix} that the diastasic potential   of $\Lambda_R$ has  the prolongation over $R$.

We consider the complexification of  $\tilde \Phi _\bold Q (z(p),\bar z(p))$ and obtain a complex holomorphic function $F_{(p,\bar p)}$ on a small neighborhood of $(p,\bar p)$ in $\Delta \subset U_X\times \bar U_X$, where $\Delta:=\{(z,\bar z)\}$.

We would like to obtain a complex-valued Hermitian positive definite function 
$\Cal F$ on    $ U_X\times U_X$
 \cite{Tr1, (2.3.2.1)} holomorphic in the first variable. If $\Cal F(z,z)$ is, in addition, positive for every $z\in U_X$ then we can define $\log \Cal F(z, z)$.
Then we will get   a positive semidefinite  Hermitian form, called the Bergman pseudo-metric
$$
\qquad \qquad ds^2_U =2 \sum g_{j k} dz_j d\overline z_k, \qquad g_{j
k}:={\partial^2\log \Cal F(z, z)\over \partial  z_j\partial\overline z_k  }.
$$ 
Recall that the positive definite $\Cal F$ is 
 a positive matrix in the sense of E.\;H.\;Moore, i.e., $\Cal F$ satisfies the property 
(ii) of \cite{Tr1, (2.3.2.1)}: 
$$
\forall\;  \tilde q_1, \dots, \tilde q_N \in U_X, \;
 \forall\;  a_1, \dots, a_N \in \bold C \implies 
\sum_{j,k}^N \Cal F(\tilde q_k,\tilde q_j)a_j \bar {a}_k \geq 0.
$$

The required $\Cal F$ will arise from the  function $F_{(p,\bar p)}$ which, in turn, arises from the diastasis of the metric  $\Lambda_\Cal L$.
First, we will establish the above property (ii) in the one-dimensional case which is the key point of the  proof. Then we will derive (ii) in general.
In short, first we get the Hilbert space $H_R$. Second, we consider the diastasis
in the one-dimensional case. Third, we consider the diastasis in the general case. Finally, we obtain the separable Hilbert space $H_{U_X}$.  

In the one-dimensional case, we get the diastasis on $R$, denoted by $D_{R}$, as well as ${1\over t}D_{R,t}$ and $\lim_{t\to \infty}{1\over t}D_{R,t}$ corresponding
$\Cal L_R^t:=\Cal L^t|_R$. 
We get  the Hilbert space $H_R$ of square-integrable holomorphic functions on $R$ and, then, the corresponding Hermitian positive definite complex-valued function on $R\times R$ (see \cite{Tr1, (2.3.2.1) and (3.2.2)}).

The general case follows from the one-dimensional case because $\pi_1(X)$ is large (compare \cite{Tr1, Appendix, (A.1.2)}). 
Suppose we are given arbitrary $N$ points $\tilde q_1, \dots, \tilde q_N$
in $U_X$ and a vector $\bold v_{\tilde{p}} \in \bold T_{\tilde{p}} U_X$, where $\tilde{p}$ is a point of $U_X$.
The diastasis on $U_X$ will produce the Hermitian positive definite function on $U_X\times U_X$ as follows. 

We consider an appropriate finite Galois  covering of $X$, denoted by $X_i\; (i\gg 0)$, containing the finite set of points $q_1, \dots, q_N\in X_i$  in a  general one-dimensional nonsingular curve (compact Riemann surface) $C_i \subset X_i$ tangent to 
$\bold v_p \in \bold T_p X_i$, where $p\in C_i$ is a given point      and $X_i$ {\it depends}\/ on $N$. 
Furthermore, we assume $\tilde q_k \mapsto q_k\; (\forall k, 1\leq  k \leq N)$, $\tilde p\mapsto p$ and $\bold v_{\tilde p} \mapsto \bold v_p$. 

We conclude the proof  as in \cite{Tr1, (4.6)}.
From the one-dimensional case, we
get the desired Hermitian positive definite function on $U_X\times U_X$.
In fact, the complex-valued {\it holomorphic}\/ function on the 
\lq\lq diagonal\rq\rq\/  $\Delta \subset U_X\times \bar{U}_X\; (\Delta:=\{(z,\bar z)\})$ obtained from the diastasic potential on $U_X$ determines a unique complex-valued function on $U_X\times U_X$ (compare \cite{Kl2, Prop.\;7.6}).

Thus we get the separable Hilbert space $H_{U_X}$, the Bergman-diastasic form $ds^2_{U_X} =2 \sum g_{j k} dz_j d\overline z_k$ and a natural  {\it immersion}\/ into a projective space.
Note that  $H_{U_X}$ is separable because the reproducing kernel determines a countable total subset of $H_{U_X}$.
The corresponding metric $ds^2_{U_X}$ arises via the immersion
$$
\Upsilon : U_X \longrightarrow \bold P(H_{U_X}^*), \quad \Upsilon^* ds^2_{\bold P(H^*_{U_X})}=ds^2_{U_X}
$$
 \cite{Kb, Chap.\;4.10, p.\;228}.
If $\{\varphi_j\}$ is an orthonormal basis of $ H_{U_X}$ then  $\Upsilon$ is
given by $u\mapsto [\varphi_0(u) : \varphi_1(u) :\dots]. $
The fundamental group $\pi_1(X)$ acts on $H_{U_X}$ as follows
$$
T_\gamma :\varphi \mapsto (\varphi\circ\gamma)\cdot Jac_\gamma \qquad (\gamma
\in \pi_1(X))
$$
where $Jac_\gamma$ is the (complex) Jacobian determinant of $T_\gamma$.
We get a $\pi_1(X)$-invariant volume form on $U_X$ defined in   
Bochner canonical coordinates as follows:
$$
\Bigl(\sum^\infty_{j=0}|\varphi_j(z)|^2 \Bigl)\!\prod^n_{\alpha=1}\!\Bigl({\sqrt{-1}\over2} dz_\alpha\!\wedge  d\bar {z}_\alpha \Bigl) = \Bigl(\sum^\infty_{j=0}|\varphi_j(z)|^2 \Bigl)\!\prod^n_{\alpha=1}\bigl(dx_\alpha\!\wedge dy_\alpha \bigl)
$$
where  $z_\alpha =x_\alpha + {\sqrt{-1}}y_\alpha$. The  Ricci form of the volume form is negative definite \cite{Kb, Chap.\;2.4.4}. It follows the canonical bundle of $X$ is ample by Kodaira.

\enddemo

In Theorem 3, we consider the following conjecture by
Kobayashi \cite{Kb}: if $X$ is a Kobayashi hyperbolic projective manifold then $\eusm K_X$ is ample. The conjecture was suggested to the author by J.-P. Demailly.

\proclaim{Theorem 3} If  $X$ is a Kobayashi hyperbolic projective manifold  then $\eusm K_X$ is ample.
\endproclaim 

\demo{Proof} We proceed by induction on $\dim X$.
Let $g_{FS}$ denote the Fubini-Study metric in  projective space.
Let $Y\subset X$ be a projective submanifold of dimension $\nu < \dim X$ with  very ample $\eusm K_Y^{t_Y}$ ($
t_Y\gg 0$). It is equipped with the real analytic Kahler metric
$$
G_Y:= \psi^*_{\eusm K^{t_Y}_Y}g_{FS} \quad  {\text{where}} \quad
 \psi_{\eusm K^{t_Y}_Y}: Y\hookrightarrow \bold P^{N_Y(t_Y)}.
$$

A section $\sigma\in H^0(\eusm K^{t_Y}_Y)$ can be expressed locally as
$\sigma=f(dy_1\wedge\cdots \wedge dy_{\nu})^{t_Y},$ where
$y_1, \dots, y_{\nu}$ is a Bochner canonical coordinate system on $Y$ with center $\bold Q$ and $f$ is a function holomorphic in the coordinate neighborhood.  Let $\sigma_0, \dots, \sigma_{N(t_Y)}$ be a basis of $H^0(\eusm K^{t_Y}_Y)$, and let 
$$
\sigma_{i_Y}=f_{i_Y}(y)(dy_1\wedge\cdots \wedge dy_\nu)^{t_Y}\quad 
(0\leq i_Y\leq N(t_Y)).
$$
Define a volume form of $G_Y$  by setting
$$
v_{G_Y}: =\biggl(\sum^{N_Y(t_Y)}_{i_Y=0} |f_{i_Y}(y)|^2\biggl)^{1/t_Y}
(\sqrt{-1})^{\nu^2} dy_1\wedge\cdots \wedge dy_{\nu}\wedge d\bar y_1
\wedge\cdots \wedge d\bar y_{\nu}.
$$

Let $D_Y(\bold Q,p)$ be the diastasic potential at $\bold Q\in Y$ (see \cite{C, Chap.\;4, especially Theorem 12}, \cite{Tr1, (2.2)}). 
It is defined (extended) everywhere except perhaps  the inverse image on $Y$ of the antipolar hyperplane of $\bold Q\in \bold P^{N_Y(t_Y)}$ (so-called infinity).
According to Calabi \cite{C, Chap.\;1, (5)}, we have the real analytic function in $p$ in a neighborhood of $\bold Q$:
$$
D_Y(\bold Q ,p)=F(y(\bold Q), \overline{y( \bold Q)})\,+\,F(y(p), \overline{y( p)})\,-\,F(y(\bold Q), \overline{y( p)})\,-\,F(y(p), \overline{y(\bold Q)}),
$$
where $ F(w, \overline{z})$ denote the  complex holomorphic function in a  neighborhood of $\bold Q\times \bar{\bold Q}$ in $Y\times \bar Y$ (compare \cite{Tr1, (2.2.1)}) arising from a potential on $Y$. Furthermore,
$$
v_{G_Y} = e^{D_Y(\bold Q,p)/t_Y}(\sqrt{-1})^{\nu^2} dy_1\wedge\cdots \wedge dy_{\nu}\wedge d\bar y_1\wedge\cdots \wedge d\bar y_{\nu}. 
$$
The associated Ricci form  Ric\,$v_{G_Y}$ is negative.

Let $\Delta$ denote the unit disk. We  consider the family
$$
 H_Y=\{(g,h)\}:={\eufm{Hol}}_{\bold Q\!\times\!\bar\bold Q}(\Delta\times\bar\Delta, 
Y\times\bar Y), \quad  g: \Delta
\rightarrow Y,\, h: \bar\Delta\rightarrow \bar Y, 0\!\times\!0\mapsto\bold Q\!\times\!\bar\bold Q,
$$
where $g$ and $h$ are holomorphic maps. We define $H_X$ by replacing $Y$ in $H_Y$ by $X$.  The families $H_Y$ and $H_X$ are equicontinuous  \cite{Kb, Theorem 2.2.23}. With a help of classical theorems of Ascoli, Arzel\'a and Montel, we will obtain a real analytic Kahler metric with the diastasic potential $D_X(\bold Q,p)$ at $\bold Q\in X$ as follows.

We consider an element $\eufm h\in H_Y$ where $\roman{im}(g)$ and $\roman{im}(h)$ 
do not intersect the infinity. Let $D_Y(\bold Q,p)$ be the diastasic potential at $\bold Q\in Y$. We get the real analytic function $D_{Y,\eufm h}(0, p)$ on $\Delta$
(by abuse of notation, now  $p\in \Delta$).

 Let  $\{Y^\alpha\}$ be a collection of all $Y$'s as above.
For the various $Y^\alpha$'s and the various elements $\eufm h_\alpha \in H_{Y^{\!\alpha}}$ as above $\eufm h$, we take the limit of functions $D_{Y^{\!\alpha},\eufm h_\alpha}(0,p)$'s  as well as $D_{Y^\alpha}(\bold Q,p)$'s.   The limit exists. 

We, then, obtain  the real analytic function denoted by $D_X(\bold Q,p)$.
 Further, if $\bold Q'\in X$ is another point then we get $D_X(\bold Q',p)$.
If $\bold Q'$ is close to $\bold Q$ (in the Kobayashi hyperbolic topology) then 
$D_X(\bold Q',p)$ is close to $D_X(\bold Q,p)$.
So, we get the real analytic Kahler metric on $X$ with the potential $D_X(\bold Q,p)$ at
$\bold Q$.

Finally, we can apply Wirtinger's theorem. Let $\omega$ be the fundamental  form
of the  Kahler metric on $X$. The induced metric on each $Y^\alpha$ 
coincides with the corresponding $G_{Y^{\!\alpha}}$.  Thus $\omega^\nu/\nu !$ restricted to $Y^\alpha$ coincides with volume form on $Y^\alpha$. 
One can compare the volume forms on $X$ and $Y^\alpha$'s as well as the  Ricci forms of the corresponding volume forms. It follows the  Ricci form of a volume form on $X$ will be negative. Hence $\eusm K_X$ is ample by Kodaira.

\enddemo

\remark{Remark 2} In Theorem 2, $U_X$ is not necessary a bounded domain in $\bold C^n$. Recently  D.\;Wu and S.-T.\;Yau \cite{ WY} have  established that a projective manifold $X$ which admits a Kahler metric with 
 negative holomorphic sectional curvature has the ample canonical 
bundle. We observe that negative holomorphic sectional curvature of a Kahler manifold
does not imply its fundamental group is nonamenable while the negative sectional  curvature yields the fundamental group is nonamenable.

\endremark

\medskip

{\it Acknowledgments.}\/ The author would like to thank Domingo Toledo for his email.
The author would like to thank Fedor Bogomolov and Fr\'ed\'eric
Campana for pointing an error  in an earlier statement of Theorem 1. The author 
had also benefited from remarks by Campana.
The last but not list, the author would like to thank Jean-Pierre Demailly
for  his encouragement.


\Refs
\widestnumber\key{Abc}

\ref
\key A \by \by  E. Abe
\book  Hopf algebras
\publ Cambridge Univ. Press, Cambridge
\yr 2004
\endref

\ref  \key  AT  \by  E. Abe, M. Takeuchi \pages  385--404
\paper Groups Associated with Some Types of Infinite Dimensional Lie Algebras
\yr1992 \vol  146
\jour   J. of Algebra. 
\endref


\ref  \key  C  \by  E. Calabi \pages  1--23
\paper Isometric imbedding of complex manifolds
\yr1953 \vol  58
\jour   Annals of Math. 
\endref

\ref 
 \key FG  \by  K. Fritzsche, H. Grauert 
\book From Holomorphic Functions to Complex manifolds
\publ Graduate Texts in Mathematics Ser., Springer 
\yr 2010
\endref




\ref
\key Kb \by S. Kobayashi 
\book  Hyperbolic Complex Spaces
\publ  Springer
\yr  1998
\endref


\ref
\key Kl1 \by J. Koll\'ar \pages 177--216 
\paper Shafarevich maps and plurigenera of algebraic varieties
\jour Invent. Math.
\yr 1993
\vol 113
\endref

\ref
\key Kl2 \bysame
\book  Shafarevich maps and automorphic forms
\publ Princeton Univ. Press, Princeton
\yr 1995
\endref

 \ref \key LS \by T. Lyons and D. Sullivan  \pages 299-323
\paper Bounded harmonic functions on coverings
\jour J. Diff. Geom.
\yr1984 \vol  19
\endref

\ref \key M \by J.\;S.\;Milne   
\book Algebraic Groups 
\yr 2015 (version 2.00)
\publ www.jmilne.org/math 
\endref

\ref \key S \by I.\;R.\;Shafarevich  \pages 214--226
\paper On some infinite-dimensional groups. $\roman{II}$ 
\jour Izv. Akad. Nauk SSSR Ser. Mat.\yr 1981  \vol 45	
\lang Russian \transl\nofrills English transl. in \jour Math. USSR-Izv. \yr 1982 \vol 18   \pages 185--194 \endref


\ref \key To1 \by D. Toledo \pages 1218-1219
\paper Bounded harmonic functions on coverings
\yr1988 \vol  104
\jour Proc. Amer. Math. Soc.
\endref

\ref \key To2 \bysame \pages 103-119
\paper Projective varieties with non-residually finite fundamental  group
\yr1993 \vol  77
\jour Inst. Hautes \'Etudes Sci. Publ. Math. 
\endref





\ref
\key Tr1 \by  R. Treger \pages
\paper Metrics on universal covering of projective variety
\jour \!arXiv:1209.3128v5.\![math.AG]\!
 \endref


\ref
\key Tr2 \bysame \pages
\paper On a conjecture of H.\;Wu
\jour arXiv:1503.00938v1.[math.AG]
 \endref

\ref\key Tr3 \bysame  \pages\paper On uniformization of compact Kahler manifolds \jour arXiv:1507.01379v3.[math.AG] \endref

\ref
\key WY \by  D. Wu and S.-T. Yau \pages 595-604
\paper Negative holomorphic curvature and positive canonical bundle
\jour Invent. Math.
\yr 2016
\vol 204
\endref


\ref \key Z \by A.\;E. Zalesskii \paper Linear groups \page 57--107 \jour Usp. Mat. Nauk
\vol {36$,$} No. 5 (221) \yr 1981 \lang Russian \transl\nofrills English transl. in 
\jour Russ. Math. Surv.  \pages 63--128\vol {36$,$} No. 5\yr 1981 \endref 
\endRefs
\enddocument